\documentclass[12pt,epsfig]{amsart}

\usepackage[margin=1.2in]{geometry}
\usepackage{amsmath,amscd}
\usepackage{amssymb}
\usepackage{amsthm}
\usepackage{amsfonts}
\usepackage{epsfig}
\usepackage{mathrsfs}
\usepackage{scrpage2}  
\usepackage{titlesec}

\titleformat{\section}{\centering\normalfont\fontsize{12}{15}\bfseries}{\thesection}{1em}{}
  
\clearscrheadfoot                 
\rehead[]{}        
\lohead[]{}        
\lehead[]{\pagemark} 
\rohead[]{\pagemark}
\cefoot[]{}  
\cofoot[]{}     

\newtheoremstyle{thmm}{1.5ex plus 1ex minus .2ex}{1.5ex plus 1ex minus
.2ex}{\rmfamily}{}{\bfseries}{}{1em}{} \theoremstyle{thmm}
\newtheorem{theorem}{Theorem}[section]

\allowdisplaybreaks

\def\refe#1{{\rm(\ref{#1})}}
\def\R{\mathbb{R}}

\def\Z{\mathbb{Z}}
\def\N{\mathbb{N}}

\def\d{\,{\rm d}}

\begin{document}

\title
{\bf A singular integral approach to the maximal 
$\bf L^p$ regularity of parabolic equations \\[25pt]}

\author{Buyang Li}



\let\thefootnote\relax\footnote{
Notes on maximal $L^p$ regularity\newline
\indent~ Department
of Mathematics, Nanjing University, Nanjing, P.R. China.
{\tt buyangli@nju.edu.cn}
}

\maketitle

\thispagestyle{empty}

%

\section{\bf Introduction}
Consider a parabolic problem
\begin{align}\label{ParabP}
\partial_tu+Au=f 
\end{align}
on a Banach space $X$.
The maximal $L^p$ regularity of the parabolic 
problem \refe{ParabP} refers to such estimates as
\begin{align}\label{MaxLp}
\|\partial_tu\|_{L^p(\R_+;X)}+\|Au\|_{L^p(\R_+;X)}
\leq C\|f\|_{L^p(\R_+;X)},\quad\mbox{for}~~1<p<\infty .
\end{align}

Traditionally, this estimate was derived for parabolic equations by using parabolic singular integrals  \cite{Calderon65, Jones64} for $X=L^p$. For example, when $Au=-\Delta u$, we have
$$
\partial_t u(t,x)=\int_0^t\int_{\R^d}k(t-s,x-y)f(s,y)\d y\d s ,
$$
where $k(t,x)=\big(\frac{|x|^2}{t^2}-\frac{d}{2t}\big)\frac{e^{-|x|^2/t}}{(2\pi t)^{d/2}}$ is a standard singular kernel of parabolic type:
$$
|\partial_t^m\partial_x^{\alpha}k(t,x)|\leq C_{\alpha,m}/(t+|x|^2)^{m+|\alpha|/2}
$$
Therefore, by the theory of singular integrals, \refe{MaxLp} holds for $X=L^p$. 

Unfortunately, this argument has not been extended to the settings of a general Banach space $X$, or more specifically, $X=L^r$ with $r\neq p$. Most existing methods in establishing \refe{MaxLp} for a general Banach space $X$ rely on the concepts of analytic semigroups \cite{CV86,HP97}, operator-valued Fourier multiplier theory and $R$-boudnedness \cite{Weis01,Weis01-2} or $H^\infty$-functional calculus \cite{Merdy98}. 

In this note, we present a simple and fundamental approach to the maximal $L^p$ regularity of parabolic problems, which only uses the concept of singular integrals of Volterra type (which are straightforward modification of the standard singular integrals).
Knowledge of analytic semigroups, $R$-boundedness or $H^\infty$-functional calculus are not required. 

\section{\bf Singular integrals of Volterra type}
Frequently we meet
the following type of operators
\begin{equation}\label{CZOdef}
Tf(t)=\int_0^t K(t,s)f(s)\d s,
\end{equation}
where $f\in L^r(\R_+\mapsto X)$, $K(t,s)$ is a map from the Banach space $X$ to the Banach space $Y$ for any fixed $0<s<t<\infty$, and the kernel $K(t,s)$ exhibits singularity such as
\begin{align}
&\|K(t,s)\|_{{\mathcal  L}(X,Y)}\leq \frac{M}{t-s},\label{khcond0}\\[5pt]
&\|K(t,s)-K(t,s_0)\|_{{\mathcal  L}(X,Y)}\leq
\frac{M|s-s_0|^\sigma}{(t-s_0)^{1+\sigma}}\quad\mbox{if}~~t-s_0\geq 2|s-s_0| , \label{khcond}\\
& \|K(t,s)-K(t_0,s)\|_{{\mathcal  L}(X,Y)}\leq
\frac{M|t-t_0|^\sigma}{(t_0-s)^{1+\sigma}}\quad\mbox{if}~~t_0-s\geq 2|t-t_0| , \label{khcond2}
\end{align}
for some $\sigma\in(0, 1]$, where $\|\cdot\|_{{\mathcal L}(X,Y)}$ denotes the operator norm. A kernel $K(t,s)$ which satisfies the
above conditions is called a standard singular kernel of Volterra type. The above
properties imply H\"{o}rmander's condition: $K(t,s)$ is integrable
in any bounded domain $U\times V$ such that $\bar U\cap \bar
V=\emptyset$ and $(t,s)\in U\times V$ implies $t>s>0$, and
\begin{align}
&\sup_{s,s_0\in\R_+}\int_{t-s_0\geq 2|s-s_0|}\|K(t,s)-K(t,s_0)\|_{{\mathcal  L}(X,Y)}\,\d t
\leq M, \label{Hormandercond}\\
&\sup_{t,t_0\in\R_+}\int_{t_0-s\geq 2|t-t_0|}\|K(t,s)-K(t_0,s)\|_{{\mathcal  L}(X,Y)}\,\d s\leq
M. \label{Hormandercond2} \medskip
\end{align}
A singular integral operator of Volterra type is a bounded linear
operator from $L^r(\R_+\mapsto X)$ to $L^r(\R_+\mapsto Y)$ for some $r\in(1,\infty]$, given by
(\ref{CZOdef}) when $f$ has compact support and $t$ does not lie in
the support of $f$, where the kernel $K(t,s)$ satisfies the
conditions (\ref{khcond0})-(\ref{khcond2}). Of course, those
conditions also imply H\"{o}rmander's conditions
(\ref{Hormandercond})-(\ref{Hormandercond2}).\medskip

\noindent{\bf Remark 2.1}\quad The above definition mimic the definition of the usual singular integral operators \cite{Grafakos08,Grafakos09,Stein70}. We should be cautious that a singular kernel of Volterra type does not define the singular integral operator directly. For example,  the kernel $K(t,s)=1/(t-s)$ satisfies  (\ref{khcond0})-(\ref{khcond2}) but the integral operator $Tf(t)=\int_0^t K(t,s)f(s)\d s$ is not defined even in the Cauchy principal sense.
\medskip

The following theorem indicates that singular integral operators of Volterra type preserve the essential properties of standard singular integral operators.

\begin{theorem}\label{MainTHM}
{\it Let $X$ and $Y$ be reflexive Banach spaces. Then a singular integral operator of Volterra type
is of
weak-type $(1,1)$ and strong-type $(p,p)$ for $1<p<\infty$. In
particular, if $\|T\|_{L^r(\R_+\mapsto X)\rightarrow L^r(\R_+\mapsto Y)}=B$, then
\begin{align*}
&\|Tf\|_{L^{1,\infty}(\R_+\mapsto Y)}\leq
C(M+B)\|f\|_{L^1(\R_+\mapsto X)},\\[3pt]
&\|Tf\|_{L^p(\R_+\mapsto Y)}\leq C_{p}(M+B)\|f\|_{L^p(\R_+\mapsto X)},\quad
1<p<\infty.
\end{align*}
}
\end{theorem}

\noindent{\it Proof}$\quad$ The theorem is proved based on the Calder\'on-Zygmund decomposition (see Appendix). 
We only prove the case
$1<r<\infty$, as the case $r=\infty$ can be proved similarly.

Without loss of generality, we can first consider $f$ as a smooth
function with compact support and then extend the result to
$L^p(\R_+\mapsto X)$ for $1\leq p<\infty$. Let $f=g+\sum_{Q_j\in \mathcal  Q}b_j$
be the Calder\'{o}n--Zygmund decomposition so that both $g$ and $b_j$
are in $L^1(\R_+\mapsto X)\cap L^\infty(\R_+\mapsto X)$ and the sum $\sum_{Q_j\in \mathcal 
Q}b_j$ converges in $L^r(\R_+\mapsto X)$. Since the operator $T$ is bounded on $L^r(\R_+\mapsto X)$, it follows that
$$Tf(t)=Tg(t)+\sum_{Q_j\in \mathcal  Q}Tb_j(t)$$ for almost all
$t\in\R_+$. The idea of such decomposition is that, if we let $Q_j^*$ be the unique cube with the same center as $Q_j$ (denoted by $s_j$), with sides parallel
to the sides of $Q_j$ and have side length $
l(Q_j^*)=2l(Q_j)$, then 
$$
Tb_j(t)=\int_0^t K(t,s)b_j(s)\d
s=\left\{
\begin{array}{ll}
\int_{Q_j}(K(t,s)-K(t,s_j))b_j(s)\d s 
    &\mbox{for $t\in (Q_j^*)^c\cap\{t>s_j\}$},\\
 0 &\mbox{for $t\in (Q_j^*)^c\cap\{t<s_j\}$} .
\end{array}
\right.
$$

Let $b=\sum_{Q_j\in{\mathcal  Q}}b_j$ and  we note that
\begin{align*}
&\big|\{t\in\R_+: \|Tf(t)\|_Y>1\}\big|\\
&\leq \big|\{t\in\R_+: \|Tg(t)\|_Y>1/2\}\big|+\big|\{t\in\R_+:
\|Tb(t)\|_Y>1/2\}\big|\\
&\leq  2^rB^r\|g\|_{L^r(\R_+\mapsto X)}^r+\big|\{t\notin\cup_j
Q_j^*: \|Tb(t)\|_Y>1/2\}\big|+\big|\cup_jQ_j^*\big|\\
&\leq
(2^{2r-1}B^r\alpha^{r-1}+C\alpha^{-1})\|f\|_{L^1(\R_+\mapsto X)}+2\sum_{j}\int_{(\cup
Q_j^*)^c}\|Tb_j(t)\|_Y\d t ,
\end{align*}
where $\alpha$ is the parameter in the Calder\'{o}n--Zygmund
decomposition. We choose $\alpha$ to satisfy
$2^{2r-1}B^r\alpha^{r-1}=2B$ so that
$$
\big|\{t\in\R_+:\|Tf(t)\|_Y>1\}\big|\leq
CB\|f\|_{L^1(\R_+\mapsto X)}+2\sum_{j}\int_{(Q_j^*)^c}\|Tb_j(t)\|_Y\d t .
$$
For the second term, we have
\begin{align*}
&\sum_{j}\int_{(Q_j^*)^c}\|Tb_j(t)\|_Y\d
t\\
&=\sum_{j}\int_{(Q_j^*)^c\cap\{t>s_j\}}\|Tb_j(t)\|_Y\d
t\\
&\leq\sum_{j}\int_{(Q_j^*)^c\cap\{t>s_j\}}\int_{Q_j}\|K(t,s)-K(t,s_j)\|_{{\mathcal  L}(X,Y)}\|b_j(s)\|_X\d
s\d t\\
&\leq \sum_{j}\sup_{s\in Q_j}\int_{t-s_j\geq
2|s-s_j|}\|K(t,s)-K(t,s_j)\|_{{\mathcal  L}(X,Y)}\d t\int_{Q_j}\|b_j(s)\|_X\d s\\
&\leq \sum_{j}M\int_{Q_j}\|b_j(s)\|_X\d s\leq CM\|f\|_{L^1(\R_+\mapsto X)}.
\end{align*}
Therefore, we have proved the weak-type $(1,1)$ estimate. The
strong-type $(p,p)$ estimates for $1<p<r$ follows from real
interpolation.

For $r<p<\infty$, we consider the transpose operator $T'$ defined by
$(Tg,f)=(g,T'f)$ for any given $g\in L^r(\R_+\mapsto X)$ and $f\in L^{r'}(\R_+\mapsto Y')$. Clearly, the operator $T'$ is given by
$$
T'f(t)=\int_t^\infty K'(t,s)f(s)\d s 
$$ with the
kernel $K'(t,s)=K(s,t)'\in {\mathcal  L}(Y',X')$. Then $T'$ is bounded from $L^{r'}(\R_+\mapsto Y')$ to
$L^{r'}(\R_+\mapsto X')$. For $f\in L^\infty(\R_+\mapsto Y')$ with compact support, we let $f=g+\sum_{Q_j\in \mathcal  Q}b_j$
be the Calder\'{o}n--Zygmund decomposition so that both $g$ and $b_j$
are in $L^1(\R_+\mapsto Y')\cap L^\infty(\R_+\mapsto Y')$ and the sum $\sum_{Q_j\in \mathcal 
Q}b_j$ converges in $L^{r'}(\R_+\mapsto Y')$, and
$$T'f(t)=T'g(t)+\sum_{Q_j\in \mathcal  Q}T'b_j(t)$$ for almost all
$t\in\R_+$, and
$$
T'b_j(t)=\int_t^\infty K(s,t)' b_j(s)\d
s=\left\{
\begin{array}{ll}
0  &\mbox{for $t\in (Q_j^*)^c\cap\{t>s_j\}$},\\
\int_{Q_j}(K(s,t)'-K(s_j,t)')b_j(s)\d s  &\mbox{for $t\in (Q_j^*)^c\cap\{t<s_j\}$} .
\end{array}
\right.
$$
Then
\begin{align*}
&\big|\{t\in\R_+: \|T'f(t)\|_{X'}>1\}\big|\\
&\leq \big|\{t\in\R_+: \|T'g(t)\|_{X'}>1/2\}\big|+\big|\{t\in\R_+:
\|T'b(t)\|_{X'}>1/2\}\big|\\
&\leq  2^rB^r\|g\|_{L^r(\R_+\mapsto Y')}^r+\big|\{t\notin\cup_j
Q_j^*: \|T'b(t)\|_{X'}>1/2\}\big|+\big|\cup_jQ_j^*\big|\\
&\leq
(2^{2r-1}B^r\alpha^{r-1}+C\alpha^{-1})\|f\|_{L^1(\R_+\mapsto Y')}+2\sum_{j}\int_{(\cup
Q_j^*)^c}\|T'b_j(t)\|_{X'}\d t ,
\end{align*}
where $\alpha$ is the parameter in the Calder\'{o}n--Zygmund
decomposition. We choose $\alpha$ to satisfy
$2^{2r-1}B^r\alpha^{r-1}=2B$ so that
$$
\big|\{t\in\R_+:\|T'f(t)\|_{X'}>1\}\big|\leq
CB\|f\|_{L^1(\R_+\mapsto Y')}+2\sum_{j}\int_{(Q_j^*)^c}\|T'b_j(t)\|_{X'}\d t .
$$
For the second term, we have
\begin{align*}
&\sum_{j}\int_{(Q_j^*)^c}\|T'b_j(t)\|_{X'}\d
t\\
&=\sum_{j}\int_{(Q_j^*)^c\cap\{t<s_j\}}\|T'b_j(t)\|_{X'}\d
t\\
&\leq\sum_{j}\int_{(Q_j^*)^c\cap\{t<s_j\}}\int_{Q_j}\|K(s,t)'-K(s_j,t)'\|_{{\mathcal  L}(Y',X')}\|b_j(s)\|_{Y'}\d
s\d t\\
&\leq \sum_{j}\sup_{s\in Q_j}\int_{s_j-t\geq
2|s-s_j|}\|K(s,t)-K(s_j,t)\|_{{\mathcal  L}(X,Y)}\d t\int_{Q_j}\|b_j(s)\|_{Y'}\d s\\
&\leq \sum_{j}M\int_{Q_j}\|b_j(s)\|_{Y'}\d s\leq CM\|f\|_{L^1(\R_+\mapsto Y')}.
\end{align*}
Therefore, we have proved the weak-type $(1,1)$ estimate. The
strong-type $(p',p')$ estimates for $1<p'<r'$ follows from real
interpolation.
In other words, $T'$ is
bounded from $L^{p'}(\R_+\mapsto Y')$ to $L^{p'}(\R_+\mapsto X')$ for
$1<p'<r'$. By a simple duality argument, this implies that $T$ is
bounded from $L^{p}(\R_+\mapsto X)$ to $L^{p}(\R_+\mapsto Y)$ for
$r<p<\infty$.

Overall, $T$ is bounded from $L^{p}(\R_+\mapsto X)$ to $L^{p}(\R_+\mapsto Y)$ for
$1<p<\infty$. \qed\medskip

\section{\bf  Maximal $\bf L^p$ regularity of parabolic equations}
\setcounter{equation}{0}

Consider the parabolic problem 
\begin{align}\label{ParabP2}
\partial_tu-\sum_{i,j=1}^d\partial_j(a_{ij}(x)\partial_i u)
+\sum_{j=1}^db_j(x)\partial_j u+c(x) u =f 
\end{align}
with the initial condition $u(0,x)\equiv 0$,
where the coefficients $a_{ij}$, $b_j$ and $c$
are bounded, measurable and satisfying the strongly ellipticity condition:
$$
\Lambda^{-1}|\xi|^2\leq \sum_{i,j=1}^da_{ij}(x)\xi_i\xi_j\leq \Lambda|\xi|^2,
\quad\forall~~x,\xi\in\R^d .
$$
This corresponds to \refe{ParabP2} with 
$Au=-\sum_{i,j=1}^d\partial_j(a_{ij}(x)\partial_i u)
+\sum_{j=1}^db_j(x)\partial_j u+c(x) u $.

Let $G(t,x,y)$ denote the Green function of the parabolic equation so that the solution of \refe{ParabP2} is given by
$$
u(t,x)=\int_0^t\int_{\R^d}G(t-s,x,y)f(s,y)\d y\d s .
$$
Thus
$$
\partial_tu(t,x)=f(t,x)+\int_0^t\int_{\R^d}\partial_tG(t-s,x,y)f(s,y)\d y\d s .
$$

Let $X=L^r$ and let the mapping from $f\in L^r(\R_+;X)$ to $\partial_tu-f\in L^r(\R_+;X)$ be denoted by $T$. 
Then $Tf(t)=\int_0^t K(t-s)f(s)\d s$ when $f$ has compact support and $t$ is not in the support of $f$, 
where the operator-valued kernel
$K(t): L^r\rightarrow L^r$ can be expressed as $$[K(t)g](x)=\int_{\R^d}\partial_tG(t,x,y)g(y)\d y.$$ 
The kernel $K(t)$ obeys
the standard estimates: 
\begin{align}
&\|K(t-s)\|_{{\mathcal  L}(L^r,L^r)}\leq \frac{C}{t-s} ,\qquad\quad\mbox{for $0<s<t<\infty$}, \label{KE1}\\
&\|\partial_tK(t-s)\|_{{\mathcal  L}(L^r,L^r)}\leq \frac{C}{(t-s)^2} , \quad\mbox{for $0<s<t<\infty$} ,
 \label{KE2}
\end{align}
which indicate that $T:L^r(\R_+;X)\rightarrow L^r(\R_+;X)$ is a singular integral operator of Volterra type. By Theorem \ref{MainTHM}, this operator must be bounded on $L^p(\R_+;X)$ for all $1<p<\infty$.
\medskip

\noindent{\bf Remark 3.1}$\quad$
By this approach, the maximal regularity 
\begin{align*}
\|\partial_tu\|_{L^p(\R_+;L^r)}+\|Au\|_{L^p(\R_+;L^r)}
\leq C\|f\|_{L^p(\R_+;L^r)},\quad\mbox{for}~~1<p,r<\infty 
\end{align*}
reduces to the homogeneous estimate
\begin{align*}
\|\partial_tu\|_{L^r(\R_+;L^r)}+\|Au\|_{L^r(\R_+;L^r)}
\leq C\|f\|_{L^r(\R_+;L^r)} ,\quad\mbox{for}~~1<r<\infty .
\end{align*}
This approach can also be applied to problems defined on a finite time interval as well as
problems defined on the bounded domain.

\section*{\bf Appendix: Calder\'{o}n--Zygmund decomposition on $\R_+$}

Let $\Z=\{0,\pm 1,\pm 2\cdots\}$ denote the set of all integers and let $\N=\{0,1,2,\cdots\}$ denote the set of all natural numbers. For any integers $n\in\Z$ and $k\in\N$ we define $Q_{n,k}=(2^nk,2^n(k+1)]$, called a dyadic cube. Then
${\mathscr Q}_+:=\{Q_{n,k}:n\in\Z~\mbox{and}~k\in\N\}$ is called the set of all dyadic cubes on the half-line $\R_+$. For any two dyadic cubes $Q_1,Q_2\in{\mathscr Q}_+$, either $Q_1\subset Q_2$ or $Q_2\subset Q_1$ or $Q_1\cap Q_2=\emptyset$. This is often referred to as the nesting property
of dyadic cubes.\\

\noindent{\bf Proposition}$\quad$ {\it Let $X$ be a
Banach space, $f\in L^1(\R_+\mapsto X)$ and
$\alpha>0$. Then there is a decomposition
$$f=g+\sum_{j=1}^\infty b_j$$
which satisfies that\\[2pt]
$(1)$ $\|g\|_{L^1(\R_+\mapsto X)}\leq \|f\|_{L^1(\R_+\mapsto X)}$ and
$\|g\|_{L^\infty(\R_+\mapsto X)}\leq 2 \alpha$,\\[2pt]
$(2)$ the functions $b_j$ are supported in disjoint dyadic cubes $Q_j\in{\mathcal  Q}_+$,
respectively,\\[2pt]
$(3)$ $\int_{Q_j}b_j(t)\d t=0$, $\int_{Q_j}\|b_j\|_X \d t\leq
4\alpha |Q_j|$,\\[2pt]
$(4)$ $\sum_j|Q_j|\leq\|f\|_{L^1(\R_+\mapsto X)} /\alpha $, \\[2pt]
$(5)$ if $f\in L^\infty(\R_+\mapsto X)$ with compact support in $\R_+$, then $g\in L^\infty(\R_+\mapsto X)$ with compact support in $\R_+$, $\int_{\R_+}\|\sum_{j=k}^lb_j\|_X^r \d t\leq 2^r\int_{\cup_{j=k}^lQ_j}\|f\|_X^r \d t$ and so the series $\sum_{j=1}^\infty b_j$ converges in $L^r(\R_+\mapsto X)$ for any $1\leq r<\infty$.\medskip
}

\noindent{\it Proof}$\quad$  Let us say that a dyadic cube $Q$ is bad if
$\frac{1}{|Q|}\int_Q\|f(t)\|_X \d t>\alpha $, and good otherwise.
A maximal bad dyadic cube is a bad dyadic cube such that any dyadic
cube strictly containing it is good. Since $f\in L^1(\R_+\mapsto X)$, any bad dyadic cube is contained in a maximal bad dyadic cube. Let ${\mathcal  Q}$ be the
collection of all maximal bad dyadic cubes. By the nesting property
of dyadic cubes, cubes in $\mathcal  Q$ are disjoint. For any $Q\in\mathcal 
Q$, $\frac{1}{|Q|}\int_Q\|f(t)\|_X \d t>\alpha $ and
$\frac{1}{|Q'|}\int_{Q'}\|f(t)\|_X \d t\leq\alpha $ for any dyadic
cube $Q'$ strictly containing $Q$. Therefore,
$$
\alpha <\frac{1}{|Q|}\int_Q\|f(t)\|_X \d t\leq 2\alpha .
$$
For any dyadic cube $Q$ outside $\cup\mathcal  Q$,
$\frac{1}{|Q|}\int_Q\|f(t)\|_X \d t\leq \alpha $. By the Lebesgue
differentiation theorem, 
$$
\|f(t)\|_X =\lim_{{\rm diam}(Q)\rightarrow
0}\frac{1}{|Q|}\int_Q\|f(s)\|_X \d s\leq 2\alpha 
$$
for almost all $t$ outside $\cup\mathcal  Q$, where $Q$ extends over all
sequence of dyadic cubes disjoint from $\cup\mathcal  Q$ and containing
$t$.

Let
$$
b_j(t)=\biggl(f(t)-\frac{1}{|Q_j|}\int_{Q_j}f(s)\d
s\biggl)1_{Q_j}(t),\qquad g(t)=\frac{1}{|Q_j|}\int_{Q_j}f(s)\d s
$$
for $t\in Q_j\subset{\mathcal  Q}$. Let $g=f$ outside $\cup\mathcal  Q$. Then
$f=g+\sum_{Q_j\in\mathcal  Q}b_j$ satisfies the requirements.
\qed


\bigskip




\begin{thebibliography}{99}

\bibitem{Calderon65}
Calder\'on, {\em Singular integrals}, Colloquium Lectures given in Aug 31--Sep 3, 1965 at the Seventieth Summer Meeting of the American Mathematical Society held in Ithaca, New York.

\bibitem{CV86}
P. Cannarsa and V. Vespri, {\em On maximal $L^p$ regularity for the abstract Cauchy problem}, Boll. Un. Mat. Ital. B, 5 (1986), pp. 165-175.

\bibitem{Grafakos08}
L. Grafakos, {\em Classical Fourier analysis}, Springer Science+Business Media, LLC, 2008.

\bibitem{Grafakos09}
L. Grafakos, {\em Modern Fourier analysis}, Springer Science+Business Media, LLC, 2009.

\bibitem{HP97}
M. Hieber and J. Pr\"{u}ss, {\em Heat kernels and maximal $L^p$-$L^q$ estimates for parabolic evolution equations}, Comm. Partial Differential Equations, 22 (1997), pp. 1647-1669.

\bibitem{Jones64}
B.F. Jones, {\em A Class of Singular Integrals},
American J. Math., 86 (1964), pp. 441-462. 

\bibitem{Merdy98}
C.L. Merdy, {\em $H^\infty$-functional calculus and applications to maximal regularity},
Publ. Math. UFR Sci. Tech. Besancon. 16 (1998), pp. 41-77.



\bibitem{Stein70}
E.M. Stein, {\em Singular integrals and differentiability properties of functions}, Princeton University Press, 1970.

\bibitem{Weis01}
L. Weis, {\em A new approach to maximal $L^p$-regularity}, Evolution equations and
their applications in physical and life sciences (Bad Herrenalb, 1998), Dekker,
New York, 2001, pp. 195-214.

\bibitem{Weis01-2}
L. Weis, {\em Operator-valued Fourier multiplier theorems and maximal $L^p$-regularity},
Math.Ann., 319 (2001), pp.735-758.



\end{thebibliography}
\end{document}